\documentclass[12pt]{article}
\usepackage{amssymb}

\begin{document}
\begin{center}
{\large\bf Attractors and repellers near generic elliptic points of reversible maps}

\vspace{12pt}

{\bf S.V.Gonchenko}$^1$, {\bf J.S.W.Lamb}$^2$, {\bf I.Rios}$^3$ and {\bf D.Turaev}$^2$\label{Author}
\vspace{6pt}

{\tiny
$^1$ Institute for Applied Mathematics \& Cybernetics, Nizhny Novgorod, Russia; E-mail: gonchenko@pochta.ru

$^2$ Mathematics Department, Imperial College, London; E-mail: jeroen.lamb@imperial.ac.uk, d.turaev@imperial.ac.uk

\vspace*{-6pt}$^3$ Instituto de Matem\'atica e Estat\'{\i}stica, Universidade Federal Fluminense; E-mail: rios@mat.uff.br}
\end{center}

{\small \em We show that resonance zones near an elliptic periodic point of a reversible map must, generically, contain
asymptotically stable and asymptotically unstable periodic orbits, along with wild hyperbolic sets.}\\

Let $M$ be a two-dimensional manifold. Let $g:M\rightarrow M$ be an
involution, so $g\circ g=id$. A map $f:M\rightarrow M$ is called {\em reversible} if it is conjugate
by $g$ to its own inverse:
\begin{equation}\label{rdf}
f^{-1}=g\circ f\circ g.
\end{equation}
Let ${\cal R}^r_g$ denote the set of the reversible maps, endowed with the $C^r$-topology, $r=1,\dots,\infty$.
In case $M$ is not compact, by $C^r$-topology we mean the topology of uniform $C^r$-convergence on compacta.
We will also consider the space ${\cal R}^\omega_{g,Q}$ of real-analytic reversible maps,
the topology on which is defined
as follows: we fix a complex neighbourhood $Q$ of $M$ and endow the space $\cal R$ by the topology of uniform convergence
on compacta in $Q$. Note that when speaking about ${\cal R}^r_g$ we assume that the involution $g$ is of class $C^r$ itself;
and when considering the space ${\cal R}^\omega_{g,Q}$ we assume $g$ is real-analytic (real on $M$ and analytic on $Q$).

A subset of $\cal R$ is called residual if it is an intersection of a countable sequence of open and dense sets; a property
is called generic if the maps for which it holds comprise a residual set.

An orbit of $f$ is called symmetric if it is invariant with respect to $g$; namely, $gx_0=f^jx_0$ for some $j$
(then, by (\ref{rdf}), $gfx_0=f^{j-1}x_0$, and so on). It is easy to see that for a symmetric
periodic orbit at least one of its points is either a fixed point of $h=g$, or a fixed point of $h=f\circ g$.
Any such point will be called a symmetric periodic point.

Dynamics in a small neighbourhood of a non-symmetric periodic orbit can be arbitrary. However, for
any such orbit $L$ its image by $g$ is also a periodic orbit and if $\lambda, \gamma$ are the multipliers of
$L$, then $\lambda^{-1},\gamma^{-1}$ are the multipliers of $Lg$. Therefore, sinks and sources in reversible
maps always exist in pairs.

For a symmetric periodic orbit, if $\lambda$ is a multiplier, then $\lambda^{-1}$ is also a multiplier.
Therefore, elliptic symmetric periodic orbits (i.e. those with multipliers $e^{\pm i\psi}$
where $\psi\neq 0,\pi$) persist at the $C^1$-small perturbations within the space of reversible maps.
It is well-known that the dynamics near a generic elliptic orbit of a reversible map appears pretty much
conservative. Namely, every point of such orbit is surrounded by invariant KAM-curves. Moreover, the
KAM-curves occupy most of the neighbourhood of the elliptic point: the relative measure
of the set filled by such curves tends to $1$ as the size of the neighbourhood tends to zero \cite{Sevr}.

However, between the KAM-curves there exist resonant zones. We show here that the generic dynamics in the
resonant zones is not conservative and have a mixed nature. Namely, our main result is the following theorem.\\

\noindent{\bf Theorem 1.} {\em In the space ${\cal R}_g^r$, $r=1,\dots,\omega$, there is
a residual subset ${\cal R}^*$ such that for every $f\in {\cal R}^*$ every point of each symmetric elliptic
periodic orbit is a limit of sinks, sources and other elliptic points.}\\

This theorem is inferred from Theorems 2 and 3 below. In order to formulate these theorems, recall the definition
of a wild hyperbolic set. Let $\Lambda$ be a non-trivial, locally-maximal, transitive, zero-dimensional, compact,
uniformly-hyperbolic set (following Newhouse, we simply call such sets basic; the restriction of the map
on a basic set is topologically conjugate to a finite Markov chain, an example is given by the Smale horseshoe).
If a smooth map has a basic set, then every $C^1$-close map has a uniquely defined basic set close to
the original one (and the two maps are topologically conjugate in a neighbourhood of the basic set).
In other words, the basic set $\Lambda$ persists at small smooth perturbations.
According to \cite{N79}, if a $C^r$-map ($r=2,...,\omega$) has a basic set $\Lambda$
such that its stable and unstable manifolds $W^s(\Lambda)$ and $W^u(\Lambda)$ {\em have a tangency} at some point (outside of
$\Lambda$), then arbitrarily $C^r$-close to the map there exists a $C^2$-open set (the Newhouse region) ${\cal N}$ such
$W^s(\Lambda)$ and $W^u(\Lambda)$, have a tangency for every map from ${\cal N}$. Following \cite{N74,N79}, we will
call the basic set $\Lambda$ {\em wild} in this case (i.e. it is wild if, both for the map itself and for every $C^2$-close
map, $W^s(\Lambda)$ and $W^u(\Lambda)$ have a tangency).

For a reversible map $f$ the basic set $\Lambda$ is called symmetric
if $g\Lambda=\Lambda$. For a symmetric wild hyperbolic set $\Lambda$ an orbit of tangency of
$W^s(\Lambda)$ and $W^u(\Lambda)$ may be symmetric with respect to $g$ (such tangencies are, generically, either quadratic
or cubic), or non-symmetric (then it is quadratic in general). The non-symmetric tangencies appear in pairs.\\~\\

\noindent{\bf Lemma 1.}{\em Let $f$ be a reversible map with a symmetric hyperbolic set $\Lambda$ whose stable and unstable manifolds
have a tangency. Then the corresponding Newhouse region $\cal N$ has a non-empty intersection ${\cal N}_{rev}$ with a neighbourhood of $f$
in the space of reversible maps. For every map from ${\cal N}_{rev}$ the manifolds $W^s(\Lambda)$ and $W^u(\Lambda)$
have a pair of non-symmetric quadratic tangencies. In ${\cal N}_{rev}$ there are dense subsets ${\cal N}_{sym,2}$
${\cal N}_{sym,3}$ such that every map from ${\cal N}_{sym,2}$ has a symmetric orbit
of quadratic tangency of $W^s(\Lambda)$ and $W^u(\Lambda)$ and every map from ${\cal N}_{sym,3}$ has a symmetric orbit
of cubic tangency of $W^s(\Lambda)$ and $W^u(\Lambda)$.}\\

The proof is obtained by noticing that near non-symmetric orbits of tangency the dynamic behaviour of a reversible map
does not differ from that for the general case, therefore the original Newhouse arguments \cite{N79} are applicable
(see also \cite{PT,GTS93b,Duarte}). The part concerning symmetric orbits is proven by showing that the symmetric tangencies
can be obtained by a perturbation of a pair of non-symmetric ones (cf. \cite{GST93a,GST99}).

Since periodic orbits are dense in $\Lambda$ and their stable and unstable manifolds are dense, respectively, in
$W^s(\Lambda)$ and $W^u(\Lambda)$, the following result holds.\\

\noindent{\bf Lemma 2.}{\em In ${\cal N}_{rev}$ the following maps form dense subsets:\\
1. reversible maps which have in $\Lambda$ a pair of non-symmetric saddle periodic orbits $P_1$ and $P_2$ such that
$P_2=gP_1$ and there is a symmetric transverse heteroclinic orbit in $W^u(P_1)\cap W^s(P_2)$ and a symmetric orbit of
quadratic heteroclinic tangency of $W^u(P_2)\cap W^s(P_1)$;\\
2. reversible maps which have in $\Lambda$ a pair of symmetric periodic orbits $Q_1$ and $Q_2$ such that
there is a non-symmetric orbit of quadratic heteroclinic tangency of $W^u(Q_1)\cap W^s(Q_2)$ and
a non-symmetric orbit of quadratic heteroclinic tangency of $W^u(Q_2)\cap W^s(Q_1)$.}\\

Note that since the orbit $P_1$ is not symmetric, one can always apply an arbitrarily small perturbation and
make the Jacobian $J(P_1)$ of the Poincare map for this
orbit (i.e. the product of the multipliers) differ from $1$ in the absolute value. Then $J(P_2)=J(P_1)^{-1}\neq 1$.
It was shown in \cite{LS} that under this condition arbitrarily small perturbations of the heteroclinic
cycle described in item 1 of Lemma 2 lead to the birth of non-symmetric periodic sinks and sources along with symmetric
elliptic periodic orbits. By \cite{DGGSL}, the same conclusion holds for perturbations of the heteroclinic cycle from item 2.
Thus, similar to the Newhouse theorem on infinitely many sinks in dissipative maps \cite{N74}, the following result
is inferred from Lemma 2 (see also \cite{GST97,LS,DGGSL}).\\

\noindent{\bf Theorem 2.} {\em In the Newhouse domain ${\cal N}_{rev}$ a residual subset $\cal B$ is comprised
by reversible maps each having infinitely many sinks, sources and symmetric elliptic points. Moreover,
a residual subset ${\cal B}^*\subset \cal B$ is comprised by the maps such that for each of them the closure of sinks,
the closure of sources and the closure of elliptic points have a non-empty intersection that contains the wild set
$\Lambda$.}\\

It is well-known that the dynamics in resonant zones near elliptic periodic points is typically chaotic. For the conservative
case the existence of transverse and non-transverse homoclinic orbits near elliptic points was shown in \cite{Z,MR,GeT}.
In the reversible case the similar statement is given by\\

\noindent{\bf Theorem 3.} {\em For a generic map from ${\cal R}_g^r$ ($r=2,\dots,\omega$) each symmetric elliptic
periodic orbit is accumulated by symmetric wild-hyperbolic sets.}\\

{\em Scheme of the proof.} It is enough to show that given any map $f\in{\cal R}_g^r$ and any elliptic periodic orbit of $f$
one can add to $f$ an arbitrarily small perturbation which does not lead it out of ${\cal R}_g^r$ and
creates a symmetric wild hyperbolic set in a given neighbourhood of this elliptic orbit.
Let $O$ be an elliptic period-$m$ point of $f$, with the multipliers $e^{\pm i\psi}$. We may always assume that $O$ is a fixed
point of the involution $h$. If $f$ is $C^r$ with $r$ finite,
we can always approximate it by a reversible $C^\infty$-map for which $O$ is a symmetric elliptic period-$m$ point,
so we need to consider further the cases $r=\infty$ or $r=\omega$ only. We can always imbed $f$ into a one-parameter family
$f_\nu$ of reversible maps for which $\psi$ changes monotonically as
the parameter varies. So, by an arbitrarily small perturbation, we may always make
$\psi(\nu_0) = 2\pi p/q$, where $p$, $q$ are mutually prime integers and $q\geq 5$. A standard fact from
normal form theory is that one can make an analytic coordinate transformation on $M$ in such a way that
the Poincare map $T= f_\nu^m$ near $O$ will, for all $\nu$ close to $\nu_0$, take the following form in the polar coordinates:
$$T = R_{2\pi p/q} \circ F + o(\rho^q)$$
where $R_\alpha$ stands for the rotation to the angle $\alpha$, and $F$ is the time-1 map of the flow
defined by
\begin{equation}\label{flowr}
\dot z=i\mu z +i \sum_{1\leq j\leq \frac{q-1}{2}} \Psi_j|z|^{2j}z+iA (z^*)^{q-1},
\end{equation}
where $z=\rho e^{i\varphi}$, $\mu=\psi(\nu)-2\pi p/q$ is a small real parameter, and $A$ and $\Psi$'s
are real. The involution $h$ with respect to which the map is reversible takes the following form
in these coordinates:
$$h: z\mapsto z^*+o(\rho^q).$$
We rewrite (\ref{flowr}) as
\begin{equation}\label{flowrp}
\begin{array}{l}
\dot \rho=A \rho^{q-1} \sin q\varphi,\\
\dot \varphi=\mu + \Psi(\rho)+A \rho^{q-2}\cos q\varphi.
\end{array}
\end{equation}
where $\Psi(\rho)=\sum_{1\leq j\leq \frac{q}{2}-1} \Psi_j \rho^{2j}$
If necessary, by an additional small perturbation we can make the coefficients $\Psi_1$ and $A$ non-zero.
By changing $\varphi\rightarrow\varphi+\pi/q$ we can always make $A\Psi_1>0$.

System (\ref{flowrp}) has $2q$ non-zero equilibria (at $\sin q\varphi=0$ and $\rho\sim\sqrt{-\mu/\Psi_1}$)
for $\mu\Psi_1<0$. They correspond to two $q$-periodic orbits of the map $T$, one of them is saddle, one is elliptic
(both the orbits are symmetric with respect to $h$).
Let $\rho^*\sim\sqrt{-\mu/\Psi_1}$ be such that $\mu=-\Psi(\rho^*)$. By scaling
$\rho=\rho^*+A(\frac{1}{2A\Psi_1}(\rho^*)^{q-2})^{1\over 2} u$ and $t=\tau/\sqrt{2A\Psi_1 (\rho^*)^q}$, we bring (\ref{flowrp})
to the form
\begin{equation}\label{penm}
\dot\varphi=u+O(|\mu|^{1/4}),\qquad \dot u=\sin q\varphi+O(|\mu|^{1/4}).
\end{equation}
The limit at $\mu\rightarrow 0$ is a pendulum equation. In this equation the separatrices of the saddle equilibria
form symmetric heteroclinic connections. As the invariant manifolds $W^u$ and $W^s$ of the saddle period-$q$ points of $T$
(in the rescaled coordinates) are close to the invariant manifolds of the saddle equilibria of (\ref{penm}) at small $\mu$,
it follows that they have symmetric homoclinic intersections (which correspond to the intersection of $W^u$ and $W^s$ with
the line of the fixed points of the involution $h$). If the symmetric homoclinic intersection stays transverse
for all small $\mu$, then the splitting angle between $W^u$ and $W^s$ must tend to zero as $\mu\rightarrow 0$, which inevitably
creates a converging to zero sequence of $\mu$ values which correspond to secondary symmetric homoclinic tangencies
(see the corresponding arguments in \cite{GeT}). Thus, there always exist arbitrarily small values of $\mu$
for which there exists a tangency (primary or secondary) between $W^u$ and $W^s$.

One can show (cf. \cite{N79,Duarte}) that arbitrarily small perturbation (within ${\cal R}_g^r$ with $r\geq 2$)
of a symmetric homoclinic tangency of invariant manifold
of a saddle periodic point creates a wild symmetric set which contains this point. This
shows that arbitrarily close to the original map $f$ there exist reversible $C^r$-maps which have
a symmetric wild hyperbolic set that contains a saddle period-$q$ point just born from the original elliptic point.

The genericity part of Theorem 3 is now obtained as follows. Take a countable base of open neighbourhoods $U_i$ in $M$.
Let ${\cal X}_i\subset{\cal R}^r_g$ be composed of such maps $F$ that either every map from some neighborhood of
$F$ in ${\cal R}^r_g$ has no elliptic points in $U_i$, or every map from some neighborhood of
$F$ in ${\cal R}^r_g$ has a symmetric wild hyperbolic set which intersects $U_i$. Each set ${\cal X}_i$ is open by construction,
moreover for every map from $\displaystyle \cap_i {\cal X}_i$ each elliptic orbit is accumulated by wild sets.
Thus, we obtain the theorem by noticing that each set ${\cal X}_i$ is dense in ${\cal R}^r_g$: if $f\not\in{\cal X}_i$,
then by an arbitrarily small perturbation an elliptic orbit can be created in $U_i$
(by the definition of ${\cal X}_i$), and then (by the arguments above)
a symmetric, persistently wild hyperbolic set can be born from this elliptic orbit. $\;\;\Box$

Theorems 2 and 3 imply Theorem 1 immediately.

In fact, Theorem 1 implies another, ``non-local" instance of coexistence of conservative
and dissipative dynamics, as given by Theorem 4 below. To formulate this theorem, let us
recall the definition of a KAM-curve. Namely, let $C$ be a closed invariant curve
of a $g$-reversible $C^r$-map $f$ ($r=4,\dots,\infty,\omega$). Assume $C$ is also invariant with respect
to the involution $g$.
Then $C$ is called a KAM-curve provided the following conditions are fulfilled.
\begin{itemize}
\item There exist $C^4$-coordinates $(\rho,\theta)$ near $C$ ($\theta$ is an angular variable and $\rho$ runs
a small interval around zero) such that $C$ is given by
$\rho=0$ in these coordinates, and the map $f$ takes the form
$$\bar \rho=q(\rho,\theta), \qquad \bar\theta=\theta+\psi(\rho,\theta),$$
with $1$-periodic in $\theta$ functions $q,\psi$ such that $q(0,\theta)\equiv 0$ and
$\psi(0,\theta)\equiv\psi_0=const.$
\item The twist condition is satisfied:
$\displaystyle\frac{\partial \psi}{\partial \rho}_{_{\rho=0}}\neq 0.$
\item The rotation number $\psi_0$ is Diophantine:
$\displaystyle \qquad|k\psi_0+p|\geq \frac{K}{|k|^\alpha}$\\
for all integer $k$ and $p$ and for some $K>0$ and $\alpha>0$ which are independent of $k,p$.
\end{itemize}

KAM-curves are robust with respect to $C^r$-small perturbations of $f$ within the class of reversible maps.
Namely, fix a sufficiently small neighbourhood $U$ of $C$. Then, for every map from ${\cal R}^r_g$
which is sufficiently close to $f$, in $U$ there exists a unique KAM-curve
with the same rotation number $\psi_0$, and this curve continuously depends on $f$.
Every KAM-curve is accumulated (from both sides)
by other KAM-curves which correspond to different values of the rotation number $\psi_0$.\\

\noindent{\bf Theorem 4.}{\em In the space ${\cal R}^r_g$ there is
a residual subset ${\cal R}^{**}$ such that for every $f\in {\cal R}^{**}$ every point of each
KAM-curve is a limit of sinks and sources.}\\

{\em Proof.} It is enough to show that given any point of any KAM-curve, one can add
an arbitrarily small perturbation after which a sink and a source is born in a given neighbourhood
of this point (the genericity part of the theorem is obtained then in the same way as it is done in Theorem 3).
As the KAM-curves do not disappear at small perturbations, we may from the very beginning assume that
the map $f$ is analytic (if it is not the case, we approximate $f$ by a $C^r$-close analytic reversible map).

Let $C_0$ be a KAM-curve with the rotation number $\psi_0$, then arbitrarily close to it there
is another KAM-curve, $C_1$, with a different rotation number, $\psi_1$, and $\psi_1-\psi_0$ is
of the order of the distance between $C_1$ and $C_0$. It is a standard fact that there must exist periodic
orbits between $C_0$ and $C_1$. Indeed, take any rational number $m/n$ between $\psi_0$ and $\psi_1$.
Then the map $F_{m/n}:(\rho,\theta)\mapsto f^n(\rho,\theta)-(0,m)$ rotates the invariant
curves $C_0$ and $C_1$ into the opposite directions, hence it must have fixed points in the annulus bounded
by these curves. By the index arguments, not all of these fixed points are saddles with positive multipliers.

Let us show that if $|n\psi_0-m|\leq 1$, then none of these fixed points
can have real negative multipliers (hence some of these fixed points are not
hyperbolic saddles). Using the Diophantine condition, one can perform averaging near
the KAM-curve, so one can show that the map $f$ near $C_0$ can be written in the form
\begin{equation}\label{kamr}
\bar\rho=\rho+O(\rho^3),\qquad \bar\theta=\theta+\Psi(\rho)+O(\rho^3),
\end{equation}
where $\Psi(0)=\psi_0$ and $\Psi'(\rho)\neq 0$. The map $F_{m/n}$ is
then $O(n\rho^2)$ close in $C^1$ to
the nonlinear rotation $(\rho,\theta)\mapsto (\rho, \theta+n\Psi(\rho)-m)$.
At a fixed point of $F_{m/n}$ we, thus, have $n\Psi(\rho)-m=O(n\rho^2)$, or
$n\psi_0-m+n\Psi'(0)\rho=nO(\rho^2)$, which implies $\rho=O(n^{-1})$. Thus, the derivative
matrix of $F_{m/n}$ at the fixed point is $O(\rho)$-close to
$\left(\begin{array}{cc} 1 & 0\\ n\Psi'(\rho) & 1\end{array}\right)$, in particular its trace
equals to $2+O(\rho)$. At small $\rho$ the trace is positive, which means that the multipliers
cannot be real negative.

Thus, we have shown that if the curves $C_0$ and $C_1$ are sufficiently close, the annulus
between $C_0$ and $C_1$ contains a fixed point of $F_{m/n}$ which is not a hyperbolic saddle.
This point is a period-$n$ point of the map $f$ and it is at a distance $O(1/n)$ from $C_0$;
moreover, by iterating (\ref{kamr}) we see that the orbit of this point gives an $O(1/n)$-approximation
of the whole of the curve $C_0$. If this non-saddle periodic orbit $L$ is not symmetric
with respect to the involution $g$, then we can always make it either a sink or a source by an arbitrarily
small perturbation of $f$. If we obtain a sink, then the image of this orbit by $g$ will be a source,
and vice versa. If the orbit $L$ is symmetric, then it is either elliptic
or parabolic. In the latter case (i.e. when it has both multipliers equal to $1$), by an arbitrarily small
perturbation of $f$ in the class of reversible maps we can make this orbit elliptic anyway.
Then, by Theorem 1, by an additional small perturbation we obtain a sink and source near this orbit.

Now, by taking $n$ sufficiently large, we obtain period-$n$ sinks and sources which approximate the KAM-curve
$C_0$ as close as we want. This proves Theorem 4.  $\;\;\Box$ \\

Note that in the analytic we do not prove here that every KAM-curve is typically accumulated by elliptic
points, this problem remains open (in the $C^\infty$-case this claim is easily proved by reducing the map to an integrable
normal form, similar to the conjecture discussed below).

We finish this paper by remark that there exists a mechanism of non-conservative behavior near
reversible elliptic points which is different from the homoclinic tangencies. Indeed, if a resonant
elliptic point is degenerate, then its perturbation leads to a pitchfork bifurcation \cite{papers,papers1,papers2}
that destroys the conservativity. Namely, like in the proof of Theorem 2, if we have a symmetric elliptic
periodic orbit with the multipliers $e^{\pm i\psi}$ where $\psi=2\pi p/q+\mu$, then for all
$\mu$ small the Poincare map near this orbit can be brought to the normal form up to the order $q+1$:
$$T = R_{2\pi p/q} \circ F + o(\rho^{q+1})$$
where $F$ is the time-1 map of the flow defined by
\begin{equation}\label{flowr1}
\dot z=i\mu z +i \sum_{1\leq j\leq \frac{q}{2}} \Psi_j|z|^{2j}z+iA (z^*)^{q-1}+iB z^{q+1}+iCz(z^*)^q,
\end{equation}
where $z=\rho e^{i\varphi}$. Note that the reversibility requires the coefficients $A$, $B$, $C$ and
$\Psi$'s be real. Importantly, for conservative maps $B$ and $C$ must be zero. However, they
can take arbitrary real values in the case of reversible systems if we do not a priori require the conservativity.
We rewrite (\ref{flowr1}) as
\begin{equation}\label{flowrp1}
\begin{array}{l}
\dot \rho=\rho^{q-1} (A+(C-B)\rho^2)  \sin \phi,\\
\dot \phi=q\mu + \sum_{1\leq j\leq \frac{q}{2}} q\Psi_j \rho^{2j}+ q\rho^{q-2}(A+(C+B)\rho^2)\cos \phi,
\end{array}
\end{equation}
where $\phi=q\varphi$. Assume the coefficient $\Psi_1$ is non-zero. If $A\neq 0$, then the behavior of system
(\ref{flowrp1}) at small $\rho$ does not differ from that of (\ref{flowrp}), i.e. it is conservative: it has only
two equilibrium states at $\mu\Psi_1<0$, and they both are symmetric (i.e. $\sin\phi=0$), so one is a saddle and one is a center.
This equilibria correspond to the pair of symmetric $q$-periodic orbits of the map $T$, one orbit is saddle and one
is elliptic. However, if $A=0$, then by changing $\mu$ and $A$ non-symmetric equilibria can be created in (\ref{flowrp1}).
Namely, if $C\neq B$, then the system has a pair of equilibria $\rho^2=A/(B-C)$ with $\sin\phi\neq 0$ for $\mu$ running a
small interval near $A\Psi_1/(C-B)$ (the boundaries of this interval correspond to pitchfork bifurcations of the symmetric
equilibria). One checks that the non-symmetric equilibria are a sink and a source if $B(B-C)>0$. They correspond
to $q$-periodic sink and source for the Poincare map $T$.

Now we make the following\\
{\bf Conjecture.} {\em Given an elliptic point with the multipliers $e^{\pm i\psi}$, an arbitrarily small (in $C^r$)
perturbation can be added such that it does not lead the map out of the class of reversible maps, and the elliptic point becomes resonant
(i.e. $\psi=2\pi p/q$) with any given real values of the coefficients $A$, $B$ and $C$ in the normal form (\ref{flowr1}).}

The problem here is to prove this conjecture for $r=\omega$, i.e. for perturbations which are small in the real-analytic sense.
For perturbations small in $C^\infty$ the conjecture is obviously true. Indeed, when $\psi/\pi$ is irrational, the normal
form for the Poincare map $T$ times the rotation $R_{-\psi}$ coincides, up to flat terms, with the time-1 shift by the orbits
of the equation $\dot z=i\Psi(zz^*)z$ for some real function $\Psi$, $\Psi(0)=0$. Therefore, by a $C^\infty$-small
perturbation of the original map
one can make the map  $R_{-\psi}\circ T$ coincide exactly with the time-$1$ shift by this flow
in a sufficiently small neighborhood of $z=0$ (though this is impossible to do by a $C^\omega$-small perturbation).
Now, take a sufficiently large $q$ such that $2\pi p/q$ is close to $\psi$, change the equation to
$\dot z=i\Psi(zz^*)z+iA (z^*)^{q-1}+iB z^{q+1}+iCz(z^*)^q$ in a sufficiently small neighbourhood of zero (so this
will be a $C^{q-2}$-small perturbation) and change the value of $\psi$ to $2\pi p/q$.

Thus, at least in the $C^\infty$-case, one can produce sinks and sources by arbitrarily small perturbation
of the elliptic points of reversible maps without the use of homoclinic tangles.

This work was supported by the RF Government grant No.11.G34.31.0039,
RFBR grants No.10-01-00429, No.11-01-00001 and No.11-01-97017-povoljie (Gonchenko),
by UK EPSRC and CARES (Brasil) grants (Lamb),
by CAPES, CNPq and Faperj (Rios), and by the Leverhulme Trust grant RPG-279 and
the RF Ministry of Education and Science grant 14.B37.21.0862 (Turaev).

\end{document}